\documentclass[11pt,a4paper,%headinclude,
fleqn,reqno]{amsart}
\usepackage[T1]{fontenc}
\usepackage[utf8]{inputenc}
\usepackage[english]{babel}
 \usepackage{csquotes}
\usepackage{graphicx}                      %
\usepackage[font=small]{quoting}
\usepackage{float}
\usepackage{caption}

\usepackage[top=1.2in,bottom=1.2in,left=1.0in,right=1.0in]{geometry}
\usepackage{verbatim}
\usepackage{color}
\usepackage{xcolor}
\usepackage{picture}
\usepackage{amsmath,amsthm,amsfonts,amssymb}
\usepackage{comment}
\usepackage{hyperref}
\hypersetup{colorlinks,linkcolor={blue}, citecolor={blue},urlcolor={black},draft}
\usepackage{enumitem}
\usepackage{mathtools}
\usepackage{wasysym}
\usepackage[url=false,doi=false,isbn=false, giveninits=true, %backref,
style=%alphabetic,
numeric,
maxbibnames=99]{biblatex}
\newtheorem{thm}{Theorem}[section]
\newtheorem{lem}[thm]{Lemma}
\newtheorem{cor}[thm]{Corollary}
\newtheorem{prop}[thm]{Proposition}

\theoremstyle{definition}

\theoremstyle{definition}

\usepackage{esint}

\usepackage{amsmath}

\def\S{\Sigma}
\def\n{\nabla}

\def\p{\partial}

\def\n{\nabla}

\def\tr{\mathrm{tr}}

\def\p{\partial}

\def\d{\delta}

\def\s{\sigma}

\def\ov{\overline}

\def\n{\nabla}
\def\<{\langle}
\def\>{\rangle}

\def\n{\nabla}

\def\RR{\mathbb{R}}

\def\SS{\mathbb{S}}

\def\tr{\mathrm{tr}}
\def\J{\mathcal{J}}

\def\p{\partial}

\def\d{\delta}

\def\s{\sigma}

\def\ov{\overline}

\def\wh{\widehat}

\def\C{\mathcal{C}}

\def\L{\mathcal{L}}

\def\wt{\widetilde}

{\left\{\begin{array}{@{}l@{}}}{\end{array}\right.}
\patchcmd{\abstract}{\scshape\abstractname}{\textbf{\abstractname}}{}{}
\makeatletter
\def\@makefnmark{}
\makeatother

\numberwithin{equation}{section}
\numberwithin{exa}{section}

\begin{document}
\title[A logarithmic Gauss curvature flow] {A capillary logarithmic Gauss curvature flow and the capillary Minkowski problem}

\author[X. Mei]{Xinqun Mei}
\address[X. Mei]{Key Laboratory of Pure and Applied Mathematics, School of Mathematical Sciences, Peking University,  Beijing, 100871, P.R. China}
\email{\href{qunmath@pku.edu.cn}{qunmath@pku.edu.cn}}

\author[L. Weng]{Liangjun Weng}
\address[L. Weng]{Centro di Ricerca Matematica Ennio De Giorgi, Scuola Normale Superiore, Pisa, 56126, Italy}
\email{\href{mailto:liangjun.weng@sns.it}  {liangjun.weng@sns.it}}



\begin{abstract}
We study a logarithmic Gauss curvature flow of strictly convex capillary hypersurfaces in the Euclidean half-space. For an acute contact angle and an arbitrary positive smooth prescribed function, we prove that a suitable initial dilation yields a global solution which converges smoothly after a horizontal translation with a uniquely determined constant velocity. As an application, we provide a logarithmic flow proof of the capillary Minkowski problem, previously established in \cite[Theorem~1.1]{MWW-AIM}.
\end{abstract}
\maketitle

%\tableofcontents

\section{Introduction}

The classical Minkowski problem asks for a convex body with prescribed surface area measure. In the smooth setting, it becomes a prescribed Gauss curvature problem and leads to a Monge--Amp\`ere equation on the unit sphere. The existence and regularity theory was established through the works of Nirenberg \cite{nirenberg1953}, Pogorelov \cite{pogorelov}, and Cheng and Yau \cite{Cheng-Yau1976}, among others. From a parabolic point of view, Chou and Wang \cite[Theorem~A]{CW00} introduced a logarithmic Gauss curvature flow for closed convex hypersurfaces. By selecting a suitable initial dilation, they obtained smooth convergence after translation and thereby gave a logarithmic flow proof of the smooth Minkowski problem. For graphs over smooth bounded strictly convex domains,
Schn\"urer and Schwetlick
\cite[Theorem~1.1]{Schnuerer-Schwetlick2004}
established long time existence and smooth convergence to
translating solutions for a non-parametric logarithmic
Gauss curvature flow with a prescribed Neumann boundary condition. 

Jointly with Wang, we introduced the capillary Minkowski problem in \cite{MWW-AIM},
which asks for a strictly convex hypersurface with prescribed Gauss
curvature and boundary meeting a supporting hyperplane at a constant
contact angle. In terms of the support function, this problem reduces
to a Monge--Amp\`ere equation on a spherical cap with a Robin boundary
condition. For $\theta\in(0, {\pi}/{2}]$, we established a necessary
and sufficient condition for solvability, together
with uniqueness up to horizontal translations
\cite[Theorem~1.1]{MWW-AIM}. The proof there uses the method of continuity and does not require any symmetry hypothesis on the prescribed function. Subsequently, related capillary Gauss curvature flows, and their volume normalized version for strictly convex capillary hypersurfaces have been studied in \cite{MWW-GCF, H-H-I}; in particular, \cite[Theorem 1.1 for  $p=1$]{H-H-I} provides a flow proof of \cite[Theorem 1.1]{MWW-AIM} under an additional capillary even symmetry hypothesis.

In this paper, we study a capillary logarithmic Gauss curvature flow for an arbitrary positive smooth prescribed function and an arbitrary smooth, strictly convex capillary initial hypersurface. For acute contact angles, we prove that a suitable initial dilation produces a global solution which converges smoothly after a horizontal translation with constant velocity. No symmetry assumption is imposed on either the prescribed function or the initial hypersurface. As a consequence, we obtain a logarithmic flow proof to recover the result of 
\cite[Theorem~1.1]{MWW-AIM}.

 Two features of the boundary problem require particular attention.
First, the principal curvature radii estimates must be compatible with the
Robin boundary condition. Second, bounds for the capillary inner and outer radii do not control the position of the evolving hypersurface in the horizontal directions. To address these issues, we combine the moving center strategy of Chou--Wang \cite{CW00} with boundary modification terms adapted to the capillary geometry. %This gives estimates independent of horizontal translations before the critical initial scale is selected. %The resulting continuation criterion and the control of the remaining translation yield convergence of the whole flow.

To state our main result, we follow the notation in \cite{MWW-AIM}
and \cite[Section~2]{mwwx}. Throughout the paper, let $n\geq2$ and
$\theta\in(0, {\pi}/{2})$. We write
\begin{eqnarray*}
 \RR^{n+1}_+
 = \left\{x\in\RR^{n+1}:x_{n+1}=\<x,E_{n+1}\>>0 \right\},
 \qquad e=-E_{n+1},
\end{eqnarray*}
where $E_i$ denotes the $i$th standard coordinate vector. A \textit{capillary hypersurface} with contact angle $\theta$ is a smooth, compact, properly embedded hypersurface $\S\subset\ov{\RR^{n+1}_+}$ satisfying $\mathrm{int}(\S)\subset\RR^{n+1}_+$, $\p\S\subset\p\RR^{n+1}_+$, and
\begin{eqnarray*} \<\nu,e\>=-\cos\theta
 \qquad\text{along }\p\S.
\end{eqnarray*}
Here $\nu$ is the outward unit normal of $\S$. For a strictly convex capillary hypersurface $\S$, we denote by $\widehat\S$ the convex body bounded by $\S$ and its planar \textit{wetting part} $\widehat{\p\S}\subset\p\RR^{n+1}_+$. For $\rho>0$ and $x_0\in\partial\ov{\RR^{n+1}_{+}}$, the spherical cap
\begin{eqnarray*}%\label{cap-definition}
 \C_{\rho,\theta}(x_0)
\coloneqq  \left\{x\in\ov{\RR^{n+1}_{+}}:
 |x-(x_0-\rho\cos\theta E_{n+1})|=\rho \right\},
\end{eqnarray*}is a capillary hypersurface with contact angle $\theta$. 
We abbreviate $\C_{\rho,\theta}(0)$ as $\C_{\rho,\theta}$ and
$\C_{1,\theta}$ as $\C_\theta$. 

Let $f\in C^\infty(\C_\theta)$ be positive.
We consider embeddings
$X(\cdot,t):M\to\ov{\RR^{n+1}_+}$
of a compact $n$-dimensional manifold $M$ with boundary.
We set $\Sigma_t\coloneqq X(M,t)$ and require
\begin{eqnarray*}
 X(\operatorname{int}M,t)\subset\RR^{n+1}_+,
 \qquad
 X(\partial M,t)\subset\partial\RR^{n+1}_+.
\end{eqnarray*}
For positive time, the embeddings are smooth and evolve by
\begin{eqnarray}\label{flow}
 \begin{cases}
 \displaystyle \p_t X
 =- \left( \log\frac{K}{f(\widetilde\nu)} \right) \widetilde\nu,
 & \qquad \text{in }M\times[0,T),\\
 \< \nu,e\>=-\cos\theta ,
 & \qquad \text{on }\partial M\times[0,T),\\
 X(\cdot,0)=X_0(\cdot),& \qquad \text{on }M.
 \end{cases}
\end{eqnarray}
Here $K$ is the Gauss curvature and $\widetilde \nu\coloneqq \nu+\cos\theta e$; see \cite[Section 2]{mwwx} for more discussion on the capillary Gauss map $\widetilde \nu$.

 We denote by $\sigma$,
$\nabla$ and $\mu$ the round metric, its Levi-Civita connection,
and the outward unit conormal of $\partial\C_\theta$, respectively.  If $\S_t$ is strictly convex, by \cite[Lemma 2.2]{mwwx}, we can parametrize $\S_{t}$ by the inverse capillary Gauss map, i.e., $X(\xi, t):\C_{\theta}\to \S_{t}$,  given by
	\begin{eqnarray*}
		X(\xi, t)= \widetilde \nu ^{-1} (\xi)  =\nu^{-1}(\xi-\cos\theta e).
	\end{eqnarray*}
The support function of $\Sigma_t$ is 
$$h(\xi,t)\coloneqq \<X(\xi,t),\xi-\cos\theta e\>. $$ In particular, the support function of $\C_\theta$ is 
\begin{eqnarray*}%\label{ell-identities}
 \ell(\xi)=\<\widetilde\nu,\nu\>
 =\sin^2\theta+\cos\theta\<\xi,e\>.
\end{eqnarray*}It is easy to check that
\begin{eqnarray}\label{eqn:ell-prop-hessian-robin}
     \n^2\ell+\ell\sigma=\sigma,
 \qquad \n_\mu\ell=\cot\theta\ell
 \quad\hbox{on }\partial\C_\theta.
\end{eqnarray}
%In particular, $\ell>0$ on $\C_\theta$.

Analogous to \cite[Proposition 3.2]{MWW-GCF}, up to tangential reparametrization,
\eqref{flow} is equivalent to
\begin{eqnarray}\label{scalarflow}
 \begin{cases}
 \partial_t h=\ell(\xi) (\log\det A+\log f)=\ell F,
 &\quad \text{in }  \C_\theta\times[0,T),\\
 \n_\mu h=\cot\theta h,
 & \quad \text{on }  \partial\C_\theta\times[0,T),\\
 h(\cdot,0)=h_0 (\cdot) ,& \quad \text{on }\C_\theta,
 \end{cases}
\end{eqnarray}
where $A\coloneqq \n^2h+h\sigma>0$, $A^{ij}$ denotes the inverse matrix
of $A$, $K=(\det A)^{-1}$ and $F\coloneqq \log\det A+\log f$.

The main result of this paper is as follows.
\begin{thm}\label{main-theorem}
Let $\theta\in(0,\pi/2 )$, let $f$ be a positive smooth
function on $\C_\theta$, and let $\S_0$ be a smooth strictly convex
capillary hypersurface.
There exists $\tau>0$ such that the solution of \eqref{flow}
starting from $\tau X_0$ exists for all time.
After reparametrization by the capillary Gauss map,
\begin{eqnarray}\label{asy}
 X(\cdot,t)+t\eta\to  X^{\ast}
 \qquad\text{ smoothly as }\quad t\to +\infty,
\end{eqnarray}
where $X^{\ast}\subset\overline{\RR^{n+1}_+}$ is a smooth strictly convex capillary hypersurface
and $\eta\in\partial {\RR^{n+1}_+}$ is uniquely determined by
\begin{eqnarray}\label{main-balance}
 \int_{\C_\theta}\frac{\xi_i}
 {e^{\ell(\xi)^{-1}\<\eta,\xi\>}f(\xi)} d\sigma=0,
 \qquad i=1,\cdots,n.
\end{eqnarray}
Moreover, the Gauss curvature of $X^{\ast}$, expressed in the capillary Gauss map 
parametrization over $\C_\theta$, is $e^{\ell(\xi)^{-1}\<\eta,\xi\>}f(\xi)$.
\end{thm}

As an application, we provide a logarithmic flow proof of the result in \cite[Theorem~1.1]{MWW-AIM}, without any additional symmetry hypothesis.

\begin{cor}
Let $\theta\in(0,\pi/2 )$. A positive smooth function $f$
on $\C_\theta$ is the Gauss curvature of a smooth strictly convex
capillary hypersurface in $\ov{\RR^{n+1}_+}$ with contact angle $\theta$ if and only if
\begin{eqnarray}\label{minkowski-balance}
 \int_{\C_\theta}\frac{\xi_i}{f(\xi)} d\sigma=0,
 \qquad i=1,\cdots,n.
\end{eqnarray}
\end{cor}

\begin{proof}
The necessity is already shown in 
\cite[Proposition~2.6]{MWW-AIM}.
If \eqref{minkowski-balance} holds, then uniqueness in Theorem~\ref{main-theorem} implies
$\eta=0$, and the limiting hypersurface has the Gauss curvature $f$.
\end{proof}

\medskip 

We briefly outline the idea of the proof for Theorem \ref{main-theorem} to conclude the introduction. Following the dilation strategy of Chou and Wang \cite{CW00}, we compare the flow with spherical caps to distinguish contracting and expanding initial data. The main estimate bounds the principal curvature radii by a subquadratic power of the outer radius. To establish this estimate, we choose a horizontal moving center and introduce a boundary modification
adapted to the Robin boundary condition, cf. Proposition \ref{prop-upper} and its proof. Combining this estimate with the capillary radius inequality in \cite{CW00,MWW-Weingarten}, we obtain uniform capillary inner and outer radii bounds for the intermediate dilations in Theorem \ref{thm-1}. A local barrier argument in the spirit of \cite{Mc2004, SW2024} and a quotient estimate of Tso's trick \cite{Tso1985} then provide an upper bound for the Gauss curvature; together, these yield the uniform parabolicity and long time existence of our flow.

To prove convergence, we determine a horizontal vector $\eta$ by the condition \eqref{main-balance}. The uniform estimates and the monotonicity formula first yield
a smooth stationary subsequential limit.
Uniqueness up to horizontal translations follows from the classical Minkowski uniqueness theorem applied to the enclosed convex bodies. Together with convergence of the weighted horizontal integrals, this identifies a single limit and yields smooth convergence of the whole translated flow. Finally, the comparison principle shows that the critical initial dilation is unique.

\bigskip

\noindent\textit{Organization of the paper.}
In Section~\ref{sec:pre}, we establish a monotonicity formula and the comparison principle, and recall the capillary radius inequality. In Section~\ref{sec:estimates}, we derive estimates for the principal curvature radii and the Gauss curvature in terms of the capillary inner and outer radii.
In Section~\ref{sec:selection}, we select a critical scaling
and establish uniform estimates.
Finally, in Section~\ref{sec:convergence}, we determine the horizontal
translation velocity, prove smooth convergence of the translated
flow, and establish uniqueness of the critical scaling.

\medskip

\noindent\textit{Acknowledgment}: The authors would like to express their sincere gratitude to Profs. X. Ma and G. Wang for their constant encouragement and support.

\medskip

\noindent\textit{Disclosure of AI use.}
By the end of 2025, the authors had completed most of the proof without using AI tools. AI tools were used only to assist in completing the proof of Proposition \ref{prop-upper}, polishing the language, and checking mathematical arguments. The authors have verified all mathematical arguments and references, reviewed all revisions, and take full responsibility for the final manuscript.

\section{Preliminaries}\label{sec:pre}

In this section, we establish two properties of the scalar flow and
recall a geometric inequality for capillary inner and outer radii. These results apply
on every interval of smooth strictly convex existence and do not
require a choice of critical initial scale. 

\subsection{A monotonicity formula}

Along our flow, we consider the functional
\begin{eqnarray}\label{energy-definition}
 \J(t)\coloneqq-{\rm Vol}(\widehat\S_t)
 +\int_{\C_\theta}\frac{h(\xi,t)}{f(\xi)} d\sigma.
\end{eqnarray}

\begin{lem}%\label{lem-mono}
Along flow \eqref{scalarflow}, the functional in
\eqref{energy-definition} is non-increasing. Moreover, 
\begin{eqnarray}\label{J}
 \frac{d}{dt}\J(t)
=-\int_{\C_\theta}\frac{\ell}{f}(e^F-1)F d\s. 
\end{eqnarray}

\end{lem}

\begin{proof}
Differentiating the Robin boundary condition in time, we obtain
\begin{eqnarray}\label{robin}
 \n_\mu\p_t h=\cot\theta \p_t h
 \quad\text{on }\partial\C_\theta.
\end{eqnarray}
By the inverse capillary Gauss parametrization,
$\nu=\xi-\cos\theta e$ is independent of time.
Differentiating $X=\n h+h\nu$, we obtain
\begin{eqnarray*}
 \p_t X=\n(\p_t h)+(\p_t h)\nu,
 \qquad
 \<\p_t X,\nu\>=\p_t h.
\end{eqnarray*}
Along $\partial\C_\theta$, we have
$$E_{n+1}=\cos\theta \nu-\sin\theta \mu.$$
Together with \eqref{robin}, we therefore find
\begin{eqnarray*}
 \<\p_t X,E_{n+1}\>
 =-\sin\theta \n_\mu\p_t h+\cos\theta \p_t h
 =0.
\end{eqnarray*}
Using the first variation
formula of volume and the area element $dS=\det A d\sigma$, 
\begin{eqnarray*}
 \frac{d}{dt}{\rm Vol}(\widehat\S_t)
=\int_{\S_t}\<\p_t X,\nu\> dS=\int_{\C_\theta}\p_t h \det A d\sigma.
\end{eqnarray*}
In view of \eqref{scalarflow} 
and $f\det A=e^F$, we obtain
\begin{eqnarray*} 
 \frac{d}{dt}\J(t)
=\int_{\C_\theta}
       \left(\frac1f-\det A\right)\p_t h d\sigma
   =-\int_{\C_\theta}\frac{\ell}{f}(e^F-1)F d\sigma
 \leq0.
\end{eqnarray*}
This completes the proof.
\end{proof}

\subsection{Comparison principle}
In the rest of the paper, we use the Einstein summation convention in local orthonormal frames;
indices after those of a tensor denote covariant derivatives.

\begin{lem}\label{lem-maxi}
Let $f_1,f_2$ be smooth positive functions on $\C_\theta$, and let
$h_1,h_2$ be smooth strictly convex solutions of flow \eqref{scalarflow}
with prescribed functions $f_1,f_2$, respectively.
If $h_1(\cdot,0)\leq h_2(\cdot,0)$ and $f_1\leq f_2$, then
$$h_1(\xi,t) \leq h_2 (\xi,t) , \qquad \forall (\xi,t)\in \C_\theta\times[0,T) .$$  
For $t>0$, the inequality is strict unless
$h_1(\cdot,0)=h_2(\cdot,0)$ and $f_1=f_2$.
\end{lem}

\begin{proof}
Let  
$\omega\coloneqq(h_1-h_2)/\ell$, and 
\begin{eqnarray*}
   A_k\coloneqq\n^2h_k+h_k\sigma, \qquad k=1,2.  
\end{eqnarray*}
In view of \eqref{scalarflow}
 and \eqref{eqn:ell-prop-hessian-robin}, we obtain
\begin{eqnarray}\label{eqn:-eq-w}
 \left\{\begin{array}{rcll}
\p_t \omega&=&a^{ij}\omega_{ij}+b^i\omega_i+c\omega
               +\log f_1-\log f_2
       &\text{in }\C_\theta\times[0,T),\\
 \n_\mu\omega&=&0
       &\text{on }\partial\C_\theta\times [0,T),\\
 \omega(\cdot,0)&=&(h_1(\cdot,0)-h_2(\cdot,0))/\ell\leq0, & \text{on } \C_\theta.
 \end{array}\right.
\end{eqnarray}
where the coefficients are 
\begin{eqnarray*}%\label{comparison-coefficients}
 a^{ij}\coloneqq \ell\int_0^1
 \big[(1-s)A_2+sA_1\big]^{ij} ds,\quad 
 b^i\coloneqq2\ell^{-1}a^{ij}\ell_j,
 \quad c\coloneqq\ell^{-1}\sum_{i=1}^n a^{ii}. 
\end{eqnarray*}
 
Fix $0<T_0<T$ and denote  $\widehat\omega\coloneqq e^{-Bt}\omega$,
where $B\geq\max\limits_{(\xi,t)\in\C_\theta\times[0,T_0]}c$.
Since $c-B\leq0$ and $f_1\leq f_2$, the maximum principle yields 
$\widehat\omega\leq0$.
If equality holds at $(\xi_0,t_0)$ with $0<t_0\leq T_0$,
the strong maximum principle and the Hopf boundary point lemma
yield $\omega\equiv0$ on $\C_\theta\times[0,t_0]$.
Thus the initial data agree, and Eq.  \eqref{eqn:-eq-w} ensures  $f_1=f_2$.
Since $T_0<T$ is arbitrary, the assertions follow.
\end{proof}

\subsection{Capillary inner
and outer radii}

Following \cite[Section~2.2]{SW2024}, for a convex capillary hypersurface $\S\subset\ov{\RR^{n+1}_+}$, we define the capillary inner and outer radii by
\begin{eqnarray*}%\label{capillary-radii}
&& r(\widehat\S,\theta)
 \coloneqq\sup\big\{\rho>0:
 \widehat{\C_{\rho,\theta}(x_0)}\subset\widehat\S
 \text{ for some }x_0\in\partial {\RR^{n+1}_+}\big\},\notag\\&&
 R(\widehat\S,\theta)
 \coloneqq\inf\big\{\rho>0:
 \widehat\S\subset\widehat{\C_{\rho,\theta}(x_0)}
 \text{ for some }x_0\in\partial {\RR^{n+1}_+}\big\}.
\end{eqnarray*}
Along the flow \eqref{flow}, we denote $r(t)\coloneqq r(\widehat\S_t,\theta)$ and
$R(t)\coloneqq R(\widehat\S_t,\theta)$. On a compact interval
$[0,T]$ where smooth existence holds, denote 
\begin{eqnarray*}%\label{radii-time}
 R_T\coloneqq\max\limits_{0\leq t\leq T}R(t),\qquad
 r_T\coloneqq\min\limits_{0\leq t\leq T}r(t).
\end{eqnarray*}
Both radii are invariant under horizontal translations.
By the capillary version of Chou--Wang's geometric inequality
\cite[Lemma~2.2]{CW00}, established in
\cite[Lemma~2.3]{MWW-Weingarten}, we have
\begin{eqnarray}\label{radius-curvature}
 \frac{R(t)^2}{r(t)}
 \leq C\max\limits_{\C_\theta}\lambda(A(\cdot,t))
 \leq C\max\limits_{\C_\theta}H(\cdot,t),
 \qquad H\coloneqq\Delta h+nh,
\end{eqnarray}
where $\lambda(A)$ is the largest eigenvalue of $A$ and $C$ depends only on $n$ and $\theta$. %Throughout the proof, $C$ denotes a positive constant which may change from line to line.  

\section{Curvature estimates}\label{sec:estimates}

In this section, we establish a priori estimates in terms of the capillary radii. We first bound the principal curvature radii by a sub-quadratic power of $R_T$, using a horizontal moving center to control translations. We then obtain an upper bound for $K$ in terms of $r_T$ and $R_T$. Together, these estimates ensure the continuation criterion used to select the initial scaling.

\subsection{Evolution equation}
 We introduce the linearized  operator
\begin{eqnarray*}%\label{linearized-operator}
 \L\coloneqq\partial_t-\ell A^{ij}\n^2_{ij}.
\end{eqnarray*}
For any $q(t)\in  \p\RR^{n+1}_+$, we define the support function with respect to $q(t)$ by
\begin{eqnarray}\label{eqn:u-centered-support-function}
u(\xi,t)\coloneqq h(\xi,t)-\<q(t),\xi\>.
\end{eqnarray}
The precise choice of $q(t)$ will be specified in the proof of Proposition \ref{prop-upper}.
 \begin{prop}
Let $h$ be a smooth strictly convex solution of flow \eqref{scalarflow}, then  
\begin{eqnarray}\label{centered-equation}
 \begin{cases}
 \p_t u=\ell F-\<q'(t),\xi\>\eqqcolon \ell F-p(\xi,t),
 &\quad \text{in }  \C_\theta\times[0,T),\\
 \n_\mu u=\cot\theta u,
 & \quad \text{on }  \partial\C_\theta\times[0,T).
 \end{cases}
\end{eqnarray}
 \end{prop}
\begin{proof}
 The first equation in \eqref{centered-equation} is straightforward, due to 
\begin{eqnarray}\label{eqn:A-hessian-u}
    A=\n^2 h+h\s= \n^2 u+u\s,~~~~\text{ in }\C_\theta. 
\end{eqnarray}On $\p\C_\theta$, we have $\mu=\sin\theta e+\cot\theta \xi$, then $\n_\mu u= \n_\mu h-\<q(t), \mu\>=\cot\theta u$. Thus the second equation in \eqref{centered-equation} follows.
\end{proof}
\begin{prop}
    Along flow \eqref{centered-equation}, there holds
\begin{flalign}
& \text{(1)}\quad
   \L u=\ell \left(F-n+u\sum_{i=1}^n A^{ii} \right)-p.
  && \label{eqn:evolution-u}\\
& \text{(2)}\quad
  \L u_k= \ell A^{ik} u_i+\ell_k F+\ell (\log f)_k-p_k,\qquad \forall~ 1\leq k\leq n. 
  && \label{eq:gradient-commutation}  
\end{flalign}
 
\end{prop}
\begin{proof}    
By \eqref{centered-equation}, a direct computation yields Eqs. \eqref{eqn:evolution-u} and \eqref{eq:gradient-commutation}.
 
\end{proof}

We next compute the evolution of $H$.
By \eqref{eqn:A-hessian-u}, we have
$H=\Delta u+nu$, the sum of the principal curvature radii of $\S_t$.
\begin{prop}%\label{trace-identities}
Along flow \eqref{scalarflow}, there holds 
\begin{eqnarray}\label{H-evo}
\L H =-\ell A^{ip}A^{jq} \n_k A_{ij}\n_k A_{pq}
       +2\<\n F,\n\ell\> -\ell H\sum_{i=1}^n A^{ii}
       +nF+\ell\left(n^2+ \Delta\log f \right).
\end{eqnarray}  

At a boundary point, choose an orthonormal frame $\{e_i\}_{i=1}^n$ with $e_n=\mu$
and diagonalizing the tangential block of $A$. Then
\begin{eqnarray}\label{H-nor}
 \n_\mu H=-\cot\theta\sum_{\alpha=1}^{n-1}
 \frac{(A_{nn}-A_{\alpha\alpha})^2}{A_{\alpha\alpha}}
 -A_{nn}\n_\mu\log f.
\end{eqnarray}
\end{prop}

\begin{proof}
By the covariant derivative commutator formula on $\C_\theta$, we have
\begin{eqnarray*}%\label{third comm}
\n_j h_{ki}=\n_kh_{ij}+h_k\delta_{ij}-h_j\delta_{ki},
\end{eqnarray*}
and \begin{eqnarray}\label{eqn:fourth-comm}
\n_{j}\n_{i}h_{pq}=\n_{q}\n_{p}h_{ij}+2(h_{pq}\delta_{ij}-h_{ij}\d_{pq})+h_{iq}\d_{pj}-h_{jp}\d_{iq}.
\end{eqnarray} 
Since $A_{ij}=h_{ij}+h\delta_{ij}$ and $H=\tr A$, \eqref{eqn:fourth-comm} implies
\begin{eqnarray}\label{h-1}
 A^{ij} \n_i\n_j H
 =\sum_{k=1}^n A^{ij}\n_{kk}^2 A_{ij}+H\sum_{i=1}^n A^{ii}-n^2.
\end{eqnarray}
Differentiating $F=\log\det A+\log f$, we get
\begin{eqnarray*}%\label{h-2}
 \Delta F=\sum_{k=1}^n
 \left(A^{ij} \n_{kk}^2 A_{ij}-A^{ip}A^{jq}\n_kA_{ij}\n_k A_{pq}\right)
 +\Delta\log f.
\end{eqnarray*}
Together with \eqref{scalarflow}, \eqref{h-1}, and using \eqref{eqn:ell-prop-hessian-robin},
\begin{eqnarray*}%\label{H-exact}
 \L H &=& \Delta(\p_t h)+n \p_t h -\ell A^{ij}\n_{ij}^2 H = \Delta(\ell F)+ n\ell F- \ell A^{ij}\n_{ij}^2 H
\notag \\&=&-\ell A^{ip}A^{jq} \n_k A_{ij}\n_k A_{pq}
       +2\<\n F,\n\ell\> -\ell H\sum_{i=1}^n A^{ii}
       +nF+n^2\ell+\ell\Delta\log f.
\end{eqnarray*}
Thus the assertion \eqref{H-evo} follows.

We proceed to show \eqref{H-nor}. By \cite[Proposition~2.8]{MWW-AIM}, on $\p\C_\theta$, there holds   
\begin{eqnarray*}%\label{h1n}
 h_{\alpha n}=A_{\alpha n}=0,
 \qquad 1\leq\alpha\leq n-1.
\end{eqnarray*}
By \cite[Eq. (3.19)]{MWW-GCF}, we have 
\begin{eqnarray}\label{Aaan}
\n_n A_{\alpha\alpha }=\cot\theta(A_{nn}-A_{\alpha\alpha}),~~1\leq \alpha\leq n-1. 
\end{eqnarray}
Taking the time derivative of Robin boundary condition in \eqref{scalarflow} yields
\begin{eqnarray*}
    \n_\mu (\ell F)= \cot\theta \ell F,
\end{eqnarray*}together with \eqref{eqn:ell-prop-hessian-robin} implies 
$0=\n_\mu F=\n_\mu (\log \det A+\log f)$, namely, 
\begin{eqnarray}\label{Annn}
 \n_n A_{nn}=-A_{nn}\left(
       \sum_{\alpha=1}^{n-1}\frac{ \n_n A_{\alpha\alpha }}{A_{\alpha\alpha}}
       +\n_\mu\log f\right).
\end{eqnarray}
In view of \eqref{Aaan} and \eqref{Annn}, we
conclude that \eqref{H-nor} holds. This completes the proof.
\end{proof}

\subsection{Principal curvature radii bound}
%\label{sec:center}

\begin{prop}\label{prop-upper}
Let $\theta\in(0,\pi/2)$, and let $f$ be a smooth positive function
on $\C_\theta$. Suppose that $h$ is a smooth strictly convex solution of \eqref{scalarflow} on $[0,T]$.
For every $\gamma\in(1,2)$, there exists a constant $C$ such that
\begin{eqnarray}\label{C2}
 H(\xi,t)\leq  C \left (1+R_T^\gamma \right),
 \qquad  \forall (\xi,t)\in\C_\theta\times[0,T].
\end{eqnarray}
Here $C$ depends only on $n$, $\theta$ $\gamma$, $f$ and
$\max\limits_{\C_\theta}H(\cdot,0)$, and is independent of $T$.
\end{prop}

\begin{proof}
To control horizontal translations, we first choose a moving center
with a speed bound depending only on an upper bound for $F$.
This construction adapts the argument in
\cite[Lemmas~2.3--2.4]{CW00} to the capillary setting.
Throughout the proof, we work on the fixed interval $[0,T]$ and denote
$$M_T\coloneqq\max\limits_{\C_\theta\times[0,T]}H.$$
The arithmetic geometric inequality yields 
\begin{eqnarray}\label{F}
 F\leq n\log(H/n)+\log f\leq C+n\log(1+H),
\end{eqnarray}
it follows 
\begin{eqnarray}\label{center-speed-input}
 F\leq C\left[1+\log(1+M_T)\right].
\end{eqnarray}
For $1\leq i\leq n$, denote the largest and smallest $i$th coordinates of the planar wetting region by
\begin{eqnarray*}
 a_i(t)\coloneqq\max_{x\in \widehat{\partial\S_t} }\<x,E_i\>,
 \qquad
 b_i(t)\coloneqq\min_{x\in \widehat{\partial\S_t} }\<x,E_i\>.
\end{eqnarray*}

Next, we express $a_i$ and $b_i$ in terms of the support function in \eqref{scalarflow}. To this end, we \textbf{claim} \begin{eqnarray}\label{eqn:extrema-support function}
 a_i(t)=\frac{h(\sin\theta E_i,t)}{\sin\theta},
 \qquad
 b_i(t)=-\frac{h(-\sin\theta E_i,t)}{\sin\theta}.
\end{eqnarray} Indeed, 
let $\bar\nu$ denote the outward unit normal of $\partial\S_{t}\subset \widehat{\partial\S_t} $
within the supporting hyperplane. Along $\partial\S_t$, the
capillary boundary condition yields 
\begin{eqnarray*}
 \nu=\sin\theta \bar\nu+\cos\theta E_{n+1},
 \qquad
 \xi=\nu+\cos\theta e=\sin\theta \bar\nu.
\end{eqnarray*}
Thus the parameter values $\xi=\pm\sin\theta E_i$ correspond
to the boundary points
\begin{eqnarray*}
 X^\pm\coloneqq X(\pm\sin\theta E_i,t),
\end{eqnarray*}
whose planar outward unit normals are $\pm E_i$, respectively.
By the convexity of $\widehat{\partial\S_t} $, for every $y\in \widehat{\partial\S_t} $, we have
\begin{eqnarray*}
 \<y-X^+,E_i\>\leq0,
 \qquad
 \<y-X^-,-E_i\>\leq0.
\end{eqnarray*}
Consequently,
$\<X^+,E_i\>=a_i(t)$ and $\<X^-,E_i\>=b_i(t)$.
Since $X^\pm$ lie in the supporting hyperplane,
$\<X^\pm,E_{n+1}\>=0$. By the definition of the support function,
we obtain
\begin{eqnarray*}
&& h(\sin\theta E_i,t)
=\<X^+,\sin\theta E_i+\cos\theta E_{n+1}\>
 =\sin\theta a_i(t),\\&& 
 h(-\sin\theta E_i,t)
=\<X^-,-\sin\theta E_i+\cos\theta E_{n+1}\>
 =-\sin\theta b_i(t).
\end{eqnarray*}
Therefore, \eqref{eqn:extrema-support function} follows. 

 \medskip 
 
Since $\ell=\sin^2\theta$ on $\partial\C_\theta$, Eqs.~\eqref{scalarflow} and \eqref{eqn:extrema-support function} imply 
\begin{eqnarray*}%\label{coordinate-speeds}
 a_i'(t)= \frac{\p_t h(\sin\theta E_i,t) }{\sin\theta}=\sin\theta F(\sin\theta E_i,t), 
\end{eqnarray*}
similarly, 
\begin{eqnarray*}
 b_i'(t)=-\sin\theta F(-\sin\theta E_i,t).
\end{eqnarray*}
In view of \eqref{center-speed-input}, we can now choose $$B_T\coloneqq C\left(1+\log(1+M_T)\right)>0 ,$$ sufficiently large such that
\begin{eqnarray}\label{eqn:bound-a-i-b-i-derivative}
    a_i' (t)\leq B_T, \qquad \text{ and }  \qquad b_i'(t)\geq-B_T, \qquad \forall 1\leq i\leq n. 
\end{eqnarray}
 
Define the horizontal point  $q(t)\coloneqq(q_1(t),\cdots, q_n(t),0)\in \p\RR^{n+1}_+$ by
\begin{eqnarray*}%\label{center-selection}
 q_i(t)\coloneqq\inf_{0\leq s\leq T}
 \left\{a_i(s)+B_T|t-s| \right\}, \qquad 1\leq i\leq n. 
\end{eqnarray*}
Then each $q_i$ is $B_T$-Lipschitz, and   
\begin{eqnarray}\label{center-box}
 b_i(t)\leq q_i(t)\leq a_i(t).
\end{eqnarray}
Indeed, taking $s=t$ proves the upper bound.
If $s\leq t$, using \eqref{eqn:bound-a-i-b-i-derivative}, then
$a_i(s)+B_T(t-s)\geq a_i(t)\geq b_i(t)$.
If $s\geq t$, using \eqref{eqn:bound-a-i-b-i-derivative} again,  then
$a_i(s)+B_T(s-t)\geq b_i(s)+B_T(s-t)\geq b_i(t)$.
Taking the infimum over $s$ proves the lower bound. Then the assertion \eqref{center-box} follows.  

By \eqref{center-box} and the diameter bound, we have
$|X-q(t)|\leq CR_T$.
After constant extension and time mollification at a sufficiently
small scale, we may assume that $q$ is smooth, preserving the
distance and speed bounds up to a fixed factor.
For this smooth choice of $q(t)$, the centered support function
$u$ defined in \eqref{eqn:u-centered-support-function} satisfies $X-q(t)=\n u+u\nu$.
By the construction of $q$ and
\eqref{eqn:bound-a-i-b-i-derivative}, we obtain
\begin{eqnarray}\label{center-estimates}
 |u|+|\n u|\leq CR_T,\qquad 
 |q'(t)| \leq C\left[1+\log(1+M_T)\right], 
\end{eqnarray} where the constant $C>0$ is independent of $T$.   
%In particular, \eqref{H-evo} and \eqref{H-nor} remain unchanged for \eqref{centered-equation}.
We point out that the estimates below do not require $u$ to be positive.

\medskip

We define $d(\xi)\coloneqq \xi_{n+1}/\sin\theta$ on $\C_\theta$.
On $\p\C_\theta$, there holds
\begin{eqnarray}\label{distance-boundary}
 d=0,\qquad \n d=-\mu.
\end{eqnarray}
Let $b>0$ be a constant to be determined later; we introduce the function
\begin{eqnarray}\label{P-radius}
 P\coloneqq H+\<\n\log f,\n u-u\n\log\ell\>
       +(1+R_T)^{\gamma-2}\Phi+b(1+R_T)d,
\end{eqnarray}
where
\begin{eqnarray*}
 \Phi\coloneqq\frac12|\n u-u\n\log\ell|^2.
\end{eqnarray*}
By \eqref{P-radius} and \eqref{center-estimates}, we obtain
\begin{eqnarray}\label{P-H}
 |P-H|\leq C(1+R_T)^\gamma.
\end{eqnarray}

On $\p\C_\theta$, choose an orthonormal frame $\{e_i\}_{i=1}^n$
with $e_n=\mu$ and diagonalizing the tangential component of $A$.
By \eqref{centered-equation} and \eqref{eqn:ell-prop-hessian-robin},
we have
\begin{eqnarray}\label{centered-boundary-1}
&& u_n-u(\log\ell)_n=0, \qquad 
 \n_n\bigl(u_\alpha-u(\log\ell)_\alpha\bigr)=0,
 \qquad 1\leq\alpha\leq n-1, 
\notag \\&& \n_n\bigl(u_n-u(\log\ell)_n\bigr)
=A_{nn}-\csc^2\theta u.
\end{eqnarray}
Consequently,
\begin{eqnarray}\label{centered-boundary}
 \n_\mu\Phi
 =\sum_{k=1}^n\bigl(u_k-u(\log\ell)_k\bigr)
       \n_n\bigl(u_k-u(\log\ell)_k\bigr)=0.
\end{eqnarray}
From \eqref{H-nor}, \eqref{distance-boundary},
\eqref{centered-boundary-1} and \eqref{centered-boundary},
using \eqref{center-estimates}, we obtain
\begin{eqnarray*}
 \n_\mu P
 &=&-\cot\theta\sum_{\alpha=1}^{n-1}
       \frac{(A_{nn}-A_{\alpha\alpha})^2}{A_{\alpha\alpha}}
       +\sum_{\alpha=1}^{n-1}(\log f)_{\alpha n}u_\alpha -\csc^2\theta u\n_\mu\log f-b(1+R_T)\\
 &\leq&CR_T-b(1+R_T)<0
 \qquad\text{on }\partial\C_\theta,
\end{eqnarray*}
where the last inequality follows by choosing $b\coloneqq C+1$,
which depends only on $n$, $\theta$ and $f$.
Thus a maximum of $P$ occurs either at $t=0$ or at a spatial interior
point $(\xi_0,t_0)$ with $t_0>0$.
If the maximum occurs at $t=0$, Eq.~\eqref{P-H} proves \eqref{C2}
and we are done. We therefore assume that $t_0>0$ in the following computation.

Now we compute the evolution equation of $P$.
Combining \eqref{eqn:evolution-u}, \eqref{eq:gradient-commutation}
and \eqref{eqn:ell-prop-hessian-robin}, we obtain
\begin{eqnarray}\label{lv-centered}
 \L\bigl(u_k-u(\log\ell)_k\bigr)
 &=&(2-\ell)A^{ik}\bigl(u_i-u(\log\ell)_i\bigr)
-2\ell(\log\ell)_kA^{ij}
       \bigl(u_i-u(\log\ell)_i\bigr)(\log\ell)_j\notag\\
 &&-u(\log\ell)_k\sum_{i=1}^nA^{ii}+n\ell(\log\ell)_k+\ell(\log f)_k-p_k+p(\log\ell)_k.
\end{eqnarray}
By the definition of $\Phi$, we have
\begin{eqnarray}\label{eqn:L-Phi-0}
 \L\Phi=\bigl(u_k-u(\log\ell)_k\bigr)
       \L\bigl(u_k-u(\log\ell)_k\bigr)-\ell A^{ij}\n_i\bigl(u_k-u(\log\ell)_k\bigr)
                  \n_j\bigl(u_k-u(\log\ell)_k\bigr).
\end{eqnarray}
By \eqref{center-estimates}, it follows from \eqref{lv-centered} that
\begin{eqnarray}\label{v-equ}
 \left|\L(\n u-u\n\log\ell)\right|
 \leq C\left(1+R_T\sum_{i=1}^nA^{ii}+|q'(t)|\right).
\end{eqnarray}
Using \eqref{eqn:A-hessian-u} and \eqref{eqn:ell-prop-hessian-robin},
a direct computation yields
\begin{eqnarray}\label{v-first}
 \n_j\bigl(u_k-u(\log\ell)_k\bigr)
 =A_{kj}-\frac{u}{\ell}\delta_{kj}
       -(\log\ell)_k\bigl(u_j-u(\log\ell)_j\bigr).
\end{eqnarray}
Together with \eqref{center-estimates}, this implies
\begin{eqnarray}\label{centered-cross}
 -\ell A^{ij}\n_i\bigl(u_k-u(\log\ell)_k\bigr)
                  \n_j\bigl(u_k-u(\log\ell)_k\bigr) &\leq& -\ell H+2nu
       +2\ell\<\n\log\ell,\n u-u\n\log\ell\>\notag\\
 &\leq&-\ell H+CR_T.
\end{eqnarray}
Inserting \eqref{v-equ} and \eqref{centered-cross}
into \eqref{eqn:L-Phi-0}, and using \eqref{center-estimates}, we obtain
\begin{eqnarray}\label{phi-evo}
 \L\Phi\leq-\ell H+CR_T\left(1+|q'(t)|\right)
       +CR_T^2\sum_{i=1}^nA^{ii}.
\end{eqnarray}
In view of \eqref{v-equ} and \eqref{v-first}, we obtain
\begin{eqnarray}\label{fv-evo}
 \left|\L\left(\<\n\log f,\n u-u\n\log\ell\>\right)\right|
 \leq C\left(1+R_T\sum_{i=1}^nA^{ii}+|q'(t)|\right).
\end{eqnarray}
By the definition of $d$, we have
\begin{eqnarray*}
 \L d=\ell\frac{\xi_{n+1}+\cos\theta}{\sin\theta}
       \sum_{i=1}^nA^{ii}.
\end{eqnarray*}
Using \eqref{F}, it follows from \eqref{H-evo} that
\begin{eqnarray}\label{H-evo-1}
 \L H\leq-\ell A^{ip}A^{jq}\n_kA_{ij}\n_kA_{pq}
       +2\<\n F,\n\ell\>-\ell H\sum_{i=1}^nA^{ii}
       +C\left[1+\log(1+H)\right].
\end{eqnarray}
To absorb the first derivative term $2\<\n F,\n\ell\>$
in \eqref{H-evo-1}, choose an orthonormal frame $\{e_i\}_{i=1}^n$
around $\xi_0$ diagonalizing $A$, with eigenvalues $\lambda_i$.
Using
$$F_k=\sum_{i=1}^n\n_kA_{ii}/\lambda_i+(\log f)_k,$$
the Cauchy--Schwarz inequality yields
\begin{eqnarray}\label{gradient-absorption}
 2\<\n F,\n\ell\>
 \leq\frac{\ell}{2}\sum_{i,k}
       {(\n_kA_{ii})^2}/{\lambda_i^2}+C
 \leq\frac{\ell}{2}A^{ip}A^{jq}
       \n_kA_{ij}\n_kA_{pq}+C.
\end{eqnarray}
At $(\xi_0,t_0)$, in view of \eqref{center-estimates},
we combine \eqref{H-evo}, \eqref{phi-evo},
\eqref{fv-evo} and \eqref{gradient-absorption} to obtain
\begin{eqnarray}\label{lp-1}
 0\leq\L P
 &=&\L H+\L\left(\<\n\log f,\n u-u\n\log\ell\>\right)
       +(1+R_T)^{\gamma-2}\L\Phi
       +b(1+R_T)\L d\notag\\
 &\leq&-\ell H\sum_{i=1}^nA^{ii}
       -\ell H(1+R_T)^{\gamma-2}
       +C_1(1+R_T)^\gamma\sum_{i=1}^nA^{ii}\notag\\
 &&+C_2(1+R_T)^{\gamma-1}
       \left[1+\log(1+M_T)\right].
\end{eqnarray}

\medskip  
By the maximality of $P$ and \eqref{P-H}, we obtain
\begin{eqnarray*} 
 M_T\leq H(\xi_0,t_0)+C_0(1+R_T)^\gamma.
\end{eqnarray*}
If $M_T\leq 2 \max \left\{C_0,\frac{C_1}{1-\cos\theta} \right\}(1+R_T)^\gamma$, then the assertion \eqref{C2} follows.
Otherwise, we have $H(\xi_0,t_0)\geq M_T/2$ and also $\ell(\xi_0) H(\xi_0,t_0)\geq C_1(1+R_T)^\gamma$,  then the third term on the right-hand side of \eqref{lp-1} is absorbed by the first term.
Consequently,
\begin{eqnarray}\label{H-log}
\frac{1-\cos\theta }2  M_T\leq C_2(1+R_T)\left[1+\log(1+M_T)\right].
\end{eqnarray}
For $M_T\geq1$, the right-hand side in \eqref{H-log} is bounded by $C(1+R_T)M_T^{(\gamma-1)/\gamma}$. It follows that $M_T\leq C(1+R_T)^\gamma$; the case $M_T<1$ is immediate. Since $(1+R_T)^\gamma\leq C(1+R_T^\gamma)$, we obtain \eqref{C2}. This completes the proof.
\end{proof}

\subsection{Upper bound for Gauss curvature}%\label{sec:gauss}

To obtain a lower bound for the principal curvature radii,
we establish an upper bound for $K$.
We combine a spherical cap barrier with Tso's quotient trick
\cite{Tso1985}.
The fixed center argument on a uniform short time interval
follows the approach of \cite{Mc2004} and its capillary
adaptation in \cite[Section~3.2]{SW2024}.
For the rest of the paper, we set
$m\coloneqq\min\limits_{\C_\theta}f$.
\begin{lem}\label{lem-contain}
Suppose that $2\rho_0<m^{-1/n}$ and 
\begin{eqnarray}\label{eqn:initial-inclusion-tso}
   \widehat{\C_{2\rho_0,\theta}(x_0)}\subset\widehat\S_{t_0}. 
\end{eqnarray} Then there exists $\delta_0>0$, depending only on $n$, $m$ and $\rho_0$, such that
\begin{eqnarray}\label{contain}
 \widehat{\C_{\rho_0,\theta}(x_0)}\subset\widehat\S_t,
 \qquad t_0\leq t<\min\{t_0+\delta_0,T\}.
\end{eqnarray}
\end{lem}
\begin{proof}
Let $\rho$ solve 
\begin{eqnarray}\label{eqn:rho-ode}\rho'(s)=\log\bigl(m\rho(s)^n\bigr),
 \qquad \rho(0)=2\rho_0.
\end{eqnarray}
Since $2\rho_0<m^{-1/n}$, the radius decreases from $2\rho_0$
to $\rho_0$ in the time
\begin{eqnarray*}
 \delta_0\coloneqq
 \int_{\rho_0}^{2\rho_0}\frac{dr}{-\log(mr^n)}>0.
\end{eqnarray*}
In particular, $\rho(s)\geq\rho_0$ for $0\leq s\leq\delta_0$,
and $\delta_0$ depends only on $n$, $m$ and $\rho_0$.

The support function of the spherical cap
$\C_{\rho(t-t_0),\theta}(x_0)$ is
\begin{eqnarray*}
 h_{\mathrm{cap}}(\xi,t)
 \coloneqq\<x_0,\xi\>+\rho(t-t_0)\ell(\xi),\qquad \xi\in \C_\theta.
\end{eqnarray*}
Since $x_0$ is horizontal, the function $\<x_0,\xi\>$
belongs to the kernel of $\n^2+\sigma$ and satisfies the
Robin boundary condition. Using
\eqref{eqn:ell-prop-hessian-robin}, we obtain
\begin{eqnarray*}
 A_{\mathrm{cap}}
 \coloneqq\n^2h_{\mathrm{cap}}+h_{\mathrm{cap}}\sigma
 =\rho(t-t_0)\sigma.
\end{eqnarray*}
Together with \eqref{eqn:rho-ode}, 
\begin{eqnarray*}
 \partial_t h_{\mathrm{cap}}
=\ell \rho'(t-t_0)=\ell\log\bigl(m\rho(t-t_0)^n\bigr)=\ell\bigl(\log\det A_{\mathrm{cap}}+\log m\bigr).
\end{eqnarray*}
Moreover, $h_{\mathrm{cap}}$ satisfies the Robin boundary condition on $\p\C_\theta$. Thus it solves \eqref{scalarflow} with prescribed function $m$.

From \eqref{eqn:initial-inclusion-tso},  
\begin{eqnarray*}
 h_{\mathrm{cap}}(\xi,t_0)
 =\<x_0,\xi\>+2\rho_0\ell(\xi)
 \leq h(\xi,t_0).
\end{eqnarray*}
Since $m\leq f$, applying Lemma~\ref{lem-maxi} from time $t_0$
yields
\begin{eqnarray*}
 h(\xi,t)
 \geq h_{\mathrm{cap}}(\xi,t) = \<x_0,\xi\>+\rho(t-t_0)\ell(\xi) \geq \<x_0,\xi\>+\rho_0\ell(\xi) , 
\end{eqnarray*}
for $t_0\leq t<\min\{t_0+\delta_0,T\}$.
Here the last inequality follows from
$\rho(t-t_0)\geq\rho_0$ and $\ell>0$.
The last expression is the support function of
$\C_{\rho_0,\theta}(x_0)$, which proves \eqref{contain}.
This completes the proof.
\end{proof}
 
We can now prove an upper bound for the Gauss curvature of our flow, which depends on the capillary inner and outer radii. 
\begin{prop}\label{Prop-K-upp}
Let $h$ be a smooth
strictly convex solution of \eqref{scalarflow} on $[0,T]$, then 
\begin{eqnarray}\label{K upper}
 K(\xi,t)\leq C, \qquad \forall (\xi,t)\in \C_\theta\times[0,T],
\end{eqnarray}
where $C$ depends on $n$, $\theta$, $f$, $r_T^{-1}$, $R_T$ and the
initial data. In particular, if the capillary inner and outer radii have
uniform positive lower and upper bounds, then $C$ is independent of $T$.
\end{prop}

\begin{proof}
Choose $\rho_0>0$ such that
$2\rho_0<\min\{r_T,m^{-1/n}\}$. For each $t_0\in[0,T)$, there is a
horizontal point $x_0$ satisfying
$\widehat{\C_{2\rho_0,\theta}(x_0)}\subset\widehat\S_{t_0}$.
After a fixed translation, we may assume that $x_0=0$ on the time
interval under consideration. By Lemma~\ref{lem-contain}, the origin
remains in the enclosed body on this interval. Together with the
diameter bound, this implies 
\begin{eqnarray}\label{support-window}
 2\varepsilon_0\leq h(\xi,t)\leq C,
 \qquad 2\varepsilon_0\coloneqq\rho_0\min\limits_{\C_\theta}\ell,
\end{eqnarray}
for $t_0\leq t<\min\{t_0+\delta_0,T\}$. Here $\rho_0$, $\varepsilon_0$
and $\delta_0$ are independent of $t_0$. Decreasing $\delta_0$ if
necessary, we also assume that $[0,\delta_0]$ lies in the initial short time existence interval.

To bound the negative part of $\p_t h $, we introduce the function 
\begin{eqnarray}\label{K-test}
Q(\xi,t)\coloneqq\frac{-(t-t_0)\p_t h (\xi,t)}{h(\xi,t)-\varepsilon_0},~~\qquad~~\forall (\xi,t)\in \C_\theta\times [t_0, \min\{t_0+\delta_0, T\}). 
\end{eqnarray}
Using Robin boundary condition for $h$ and \eqref{robin}, we obtain
\begin{eqnarray*}
 \n_\mu Q=-\frac{\varepsilon_0\cot\theta}{h-\varepsilon_0}Q
 \qquad\hbox{on }\partial\C_\theta.
\end{eqnarray*}
Since $\theta\in(0,\pi/2)$, a positive maximum value of $Q$ cannot occur on the boundary. Also, $Q(\cdot,t_0)=0$. Thus a positive maximum on
any compact subinterval is attained at an interior point with
$s\coloneqq t-t_0>0$.

Differentiating \eqref{scalarflow} in time yields
$$\L (\p_t h) =\ell \left(\sum_{i=1}^n A^{ii} \right)\p_t h.$$ At the interior maximum point, $\n Q=0$, therefore yielding
\begin{eqnarray}\label{LP}
 0 \leq \L Q &= &- \frac{\p_t h}{h-\varepsilon_0}- \frac{s}{h-\varepsilon_0}\L (\p_t h) +s\frac{\p_t h}{(h-\varepsilon_0)^2} \L h   +\frac{2\ell}{h-\varepsilon_0} A^{ij}\n_i Q \n_j h
\notag  \\&=& \frac{Q+Q^2}{s}
 +\frac{n\ell-\varepsilon_0\ell\sum_{i=1}^n A^{ii}}{h-\varepsilon_0}Q.
\end{eqnarray}
On the other hand, the arithmetic geometric inequality and
\eqref{K-test} imply 
\begin{eqnarray}\label{AG}
 \sum_{i=1}^n A^{ii} \geq n(\det A)^{-1/n} = nf^{1/n}\exp\left(\frac{(h-\varepsilon_0)Q}{n\ell s}\right).
\end{eqnarray}
Dividing \eqref{LP} by $Q/s$ and inserting \eqref{support-window} and
\eqref{AG}, we obtain
\begin{eqnarray}\label{K-test-estimate}
 0\leq 1+Q+Cs-C^{-1}s\exp(C^{-1}Q/s).
\end{eqnarray}
For $C^{-1}Q\geq\delta_0$, the function
$s\mapsto s\exp(C^{-1}Q/s)$ is decreasing on $(0,\delta_0]$.
Hence the last term in \eqref{K-test-estimate} dominates the other
terms as $Q\to+ \infty$, uniformly for $0<s\leq\delta_0$. We conclude
that $Q\leq C$, with $C$ independent of $t_0$. For $t_0+\delta_0/2\leq t<\min\{t_0+\delta_0,T\}$,
\eqref{K-test} and \eqref{support-window} imply
\begin{eqnarray*}
 -\p_t h \leq \frac{C(h-\varepsilon_0)}{t-t_0}
 \leq\frac{2C(h-\varepsilon_0)}{\delta_0}
 \leq C.
\end{eqnarray*}
Consequently,
\begin{eqnarray*}
 K=(\det A)^{-1}
 =f\exp(-\ell^{-1}\p_t h)\leq C
 \qquad\text{for }\
 t_0+\delta_0/2\leq t<\min\{t_0+\delta_0,T\}.
\end{eqnarray*}

 For each $t\geq\delta_0/2$, we apply
this estimate with $t_0=t-\delta_0/2$. The remaining initial interval
is controlled by short-time existence and the initial hypersurface.
The bound at $t=T$ follows by continuity.
This proves \eqref{K upper}.
\end{proof}

\subsection{Regularity estimate}\label{sec:regularity}

For strictly convex initial support functions satisfying the Robin
boundary condition, short-time existence and uniqueness follow from
the standard theory; see \cite{Dong,L96}. Suppose that $R_T$ and $r_T^{-1}$ are bounded. Combining
Propositions~\ref{prop-upper} and \ref{Prop-K-upp}, we obtain upper
bounds for the eigenvalues of $A$ and a positive lower bound for
$\det A$. Together with \eqref{scalarflow}, this yields 
\begin{eqnarray}\label{uniform-parabolicity}
 C^{-1}\sigma\leq A=\n^2 h+h\sigma\leq C\sigma,
 \qquad |\p_t h|\leq C .%, \qquad \frac{\partial(\ell F)}{\partial h_{ij}}=\ell A^{ij}.
\end{eqnarray}
Thus the equation is uniformly parabolic and concave. The Robin boundary condition is uniformly oblique and
equivalently reads
\begin{eqnarray}\label{neumann-regularity}
 \n_\mu(\ell^{-1}h)=0
 \qquad\hbox{on }\partial\C_\theta.
\end{eqnarray}

To control the horizontal velocity, fix $I=[t_0,t_0+2]$ within the
existence interval and choose a fixed point $x_0$ in the planar wetting part
at time $t_0$. The diameter bound at $t_0$ and the bound for $\p_t h$
yield
\begin{eqnarray}\label{fixed-center-bound}
 \sup_{(\xi,t)\in \C_\theta\times I}|h(\xi,t)-\<x_0,\xi\>|\leq C.
\end{eqnarray}
Since horizontal linear functions lie in the kernel of
$\n^2+\sigma$, \eqref{uniform-parabolicity} also controls the Hessian
of $h-\<x_0,\xi\>$, and interpolation bounds its gradient.
The same argument applies to intervals truncated at a finite endpoint.

Define $\wh h \coloneqq\ell^{-1}(h-\<x_0,\xi\>)$. \eqref{neumann-regularity} implies  
\begin{eqnarray}\label{eqn:wh-h-bdry}
    \n_\mu \wh h=0,\qquad \text{ on } \p\C_\theta. 
\end{eqnarray}The preceding estimates \eqref{uniform-parabolicity} and \eqref{fixed-center-bound} imply  
\begin{eqnarray}\label{eqn:bound-u-c2-time}
\sup_{t\in I}
\left(\|\wh h(\cdot,t)\|_{C^2}
+\|\partial_t \wh h(\cdot,t)\|_{C^0}\right)\leq C.
\end{eqnarray}
For $s,t\in I$, we have 
$\|\wh h(\cdot,t)-\wh h(\cdot,s)\|_{C^0(\C_\theta)}\leq C|t-s|$ and
$\|\wh h(\cdot,t)-\wh h(\cdot,s)\|_{C^2(\C_\theta)}\leq C$.
By spatial interpolation on the fixed smooth cap, we obtain
\begin{eqnarray}\label{gradient-time-holder}
 \|\n \wh h(\cdot,t)-\n \wh h(\cdot,s)\|_{C^0(\C_\theta)}
 &=&\|\n (\wh h(\cdot,t)-\wh h(\cdot,s))\|_{C^0(\C_\theta)}\notag\\
 &\leq&C\|\wh h(\cdot,t)-\wh h(\cdot,s)\|_{C^0(\C_\theta)}^{1/2}
          \|\wh h(\cdot,t)-\wh h(\cdot,s)\|_{C^2(\C_\theta)}^{1/2}\notag\\
 &\leq&C|t-s|^{1/2},\qquad \forall s,t\in I. 
\end{eqnarray}
Together with the spatial $C^2$ and time Lipschitz bounds, this yields
uniform $C^{1+\alpha,(1+\alpha)/2}$ estimates for $\wh h$ and
$h-\<x_0,\xi\>$ for every fixed $\alpha\in(0,1)$.

 We now apply boundary regularity to the equation for $\wh h$.
Using \eqref{eqn:ell-prop-hessian-robin}, we obtain
\begin{eqnarray*}
 \partial_t\wh h=\log\det\left(
 \ell\n_{ij}^2\wh h+\ell_i\n_j\wh h+\ell_j\n_i\wh h
 +\wh h\sigma_{ij}\right)+\log f,
\end{eqnarray*}
with \eqref{eqn:wh-h-bdry}. 
By \eqref{eqn:bound-u-c2-time} and \eqref{gradient-time-holder},
the lower order arguments $\wh h$ and $\n\wh h$ have uniform
parabolic H\"older bounds. On the matrix range in
\eqref{uniform-parabolicity}, the operator is uniformly elliptic
and concave. After a uniformly elliptic concave extension outside
this matrix range, we may apply the localized boundary second
derivative H\"older estimate in local boundary coordinates;
see \cite[Theorem~14.22]{L96} and \cite[Theorem~5]{Ura}.
We obtain a uniform $C^{2+\beta,1+\beta/2}$ bound for $\wh h$,
and hence for $h-\<x_0,\xi\>$, on the latter half of $I$,
for some $\beta\in(0,1)$.
Linear parabolic Schauder estimates and differentiation then
yield all higher order bounds on smaller time intervals away
from the initial time of $I$. These bounds are independent
of $t_0$ and of the chosen horizontal center.

Consequently, a solution extends beyond a finite maximal time if
$R_T$ and $r_T^{-1}$ remain bounded. In this finite time case,
\eqref{uniform-parabolicity} also bounds the displacement from any
fixed initial horizontal center. If the radius bounds hold for all
time, the same regularity estimates are uniform on the latter halves
of the fixed length intervals, after the translations used above.

\section{The critical scaling of initial data}\label{sec:selection}

Using the preceding estimates, we now select an initial scale for which the solution remains uniformly bounded in radius. We follow the threshold argument of Chou--Wang \cite{CW00}, with spherical caps as comparison solutions in our capillary setting.

After a fixed horizontal translation, we may assume that the origin
lies in the relative interior of the planar wetting part of
$\wh\Sigma_0$. Then $h_0>0$ on $\C_\theta$.
For each dilation, undoing this normalization changes the solution
only by a fixed horizontal translation.
The capillary inner and outer radii of the spherical cap barriers solve
\begin{eqnarray}\label{cap-barrier-odes}
 \rho'(t)=\log\left(\rho(t)^n \min\limits_{\C_\theta}f\right)
 \quad\hbox{or}\qquad
\rho'(t)=\log\left(\rho(t)^n\max\limits_{\C_\theta}f\right).
\end{eqnarray}
By Lemma~\ref{lem-maxi}, the first equation provides inner barriers
and the second provides outer barriers. A subcritical outer barrier,
whose initial radius is less than $(\max\limits_{\C_\theta}f)^{-1/n}$,
shrinks to its horizontal center in finite time. A supercritical inner
barrier, whose initial radius exceeds $(\min\limits_{\C_\theta}f)^{-1/n}$, exists for all time
and expands without bound. Supercritical outer barriers also exist
for all time and have finite radius on every finite interval. These
properties follow by separating variables in \eqref{cap-barrier-odes}.

\begin{lem}\label{finite time-extinction}
If $T<\infty$ is the maximal existence time of a strictly convex
solution to \eqref{scalarflow}, then $R(t)\to0$ as $t\to T$. Moreover, the enclosed bodies $\wh\S_t$ converge to a single point in the supporting hyperplane $\p\RR^{n+1}_+$ as $t\to T$.
\end{lem}

\begin{proof}
Enclose the initial body $\wh\S_0$ in a spherical cap with radius greater than
$(\max\limits_{\C_\theta}f)^{-1/n}$. Comparison with the corresponding
expanding outer barrier confines the solution to a bounded region
on $[0,T)$ and yields
\begin{eqnarray*}%\label{finite time-radius-upper}
 \sup_{0\leq t<T}R(t)<\infty.
\end{eqnarray*}
Together with \eqref{C2}, we know $H$ is uniformly bounded on this
interval. We next \textbf{claim} that
\begin{eqnarray}\label{radius-liminf-zero}
 \liminf_{t\to T}R(t)=0.
\end{eqnarray}
We argue by contradiction. If \eqref{radius-liminf-zero} does not hold, then positivity and continuity
on compact subintervals provide a positive lower bound for $R(t)$ on
$[0,T)$. Inserting this bound into \eqref{radius-curvature}, we obtain
a positive lower bound for $r(t)$. The regularity estimates in
Section~\ref{sec:regularity} then extend the solution beyond $T$,
which is a contradiction. Hence the assertion 
\eqref{radius-liminf-zero} follows. 

Choose $t_j\to T$ such that $R(t_j)\to0$.
For sufficiently large $j$, the convex body $\wh\S_{t_j}$ is contained in a spherical cap of radius
$2R(t_j)<(\max\limits_{\C_\theta}f)^{-1/n}$. Compare the solution with the
shrinking outer barrier starting from this cap. The solution cannot
remain smooth and strictly convex past the extinction time of that
barrier. Consequently, comparison holds throughout its remaining
existence interval and yields
\begin{eqnarray*}
 R(t)\leq2R(t_j),\qquad t_j\leq t<T.
\end{eqnarray*}
Letting $j\to + \infty$, we conclude that $R(t)\to0$.

\

Let $\rho_j\leq CR(t_j)$ be the initial radius of the $j$th shrinking
spherical-cap barrier, and let $x_j$ be its fixed horizontal center.
By the comparison principle,
\begin{eqnarray*}
 \widehat \S_t\subset\wh{\C_{\rho_j,\theta}(x_j)},
 \qquad t_j\leq t<T.
\end{eqnarray*}
For any $X_{0}\in\wh{\C_{\rho_j,\theta}(x_j)}$, the definition of the cap implies
\begin{eqnarray*}
 |X_{0}-x_j|^2
 \leq \rho_j^2\sin^2\theta
       -2\rho_j\cos\theta\<X_{0},E_{n+1}\> \leq \rho_j^2\sin^2\theta.
\end{eqnarray*}
For $\max\{t_j,t_k\}\leq t<T$, choose $X\in\widehat \S_t$. Then
\begin{eqnarray*}
 |x_j-x_k|\leq|x_j-X|+|X-x_k|\leq(\rho_j+\rho_k)\sin\theta.
\end{eqnarray*}
Since $\rho_j\to0$, the sequence $\{x_j\}$ converges to a point
$x_\infty$ in the supporting hyperplane. Letting $k\to\infty$, we obtain
\begin{eqnarray*}
 |x_j-x_\infty|\leq\rho_j\sin\theta.
\end{eqnarray*}

To compare the enclosing spherical caps with a fixed horizontal center, we use
the horizontality of $x_j-x_\infty$ to derive
\begin{eqnarray*}
 |x_j-x_\infty+\rho_j\cos\theta E_{n+1}|^2
=|x_j-x_\infty|^2+\rho_j^2\cos^2\theta\leq\rho_j^2.
\end{eqnarray*}
Thus the distance between the centers of the full balls defining
$\wh{\C_{\rho_j,\theta}(x_j)}$ and
$\wh{\C_{2\rho_j,\theta}(x_\infty)}$ is at most $\rho_j$.
The corresponding ball inclusion therefore yields
\begin{eqnarray*}
 \widehat \S_t
 \subset\wh{\C_{\rho_j,\theta}(x_j)}
 \subset\wh{\C_{2\rho_j,\theta}(x_\infty)},
 \qquad t_j\leq t<T.
\end{eqnarray*}
For every $\varepsilon>0$, choose $j$ sufficiently large that
$2CR(t_j)<\varepsilon$. Since spherical caps with the same horizontal center
are nested with respect to their radii, we obtain
\begin{eqnarray*}
 \widehat \S_t\subset\wh{\C_{\varepsilon,\theta}(x_\infty)},
 \qquad t_j\leq t<T.
\end{eqnarray*}
Moreover,
$\sup\limits_{X\in\widehat \S_t}|X-x_\infty|
\leq2\rho_j\sin\theta\to0$.
Hence $\widehat \S_t$ converges to $\{x_\infty\}$ in Hausdorff
distance as $t\uparrow T$. This completes the proof.

\end{proof}

Combining this lemma with spherical cap comparison, we obtain two
criteria that can be applied at any finite time:
\begin{eqnarray}\label{radius-triggers}
 R(t_0)<\left(\max\limits_{\C_\theta}f\right)^{-1/n}
 &\Longrightarrow&\text{finite time extinction at a point},\notag\\
 r(t_0)>m^{-1/n}
 &\Longrightarrow&\text{global existence and }r(t)\to\infty.
\end{eqnarray}
For the first implication, enclose the body in a strictly subcritical
outer cap and apply Lemma~\ref{finite time-extinction}. For the second,
insert a strictly supercritical inner cap. Its positive and increasing
radius rules out the finite time extinction described in that lemma.
The solution is therefore global, and the inner comparison implies
$r(t)\to\infty$.

For the solution with initial support function $\tau h_0$, define
\begin{eqnarray}\label{critical-scales}
&&  \tau_{\ast}
 \coloneqq \sup\{\tau>0:\text{the solution shrinks to a point in finite time}\},
 \notag\\&& 
 \tau^{\ast}
  \coloneqq \inf\{\tau>0:\text{the solution is global and }r(t)\to\infty\}.
\end{eqnarray}
Both sets are nonempty. Indeed, comparison of $\tau h_0$ with
multiples of $\ell$ makes the first criterion in
\eqref{radius-triggers} applicable for small $\tau$ and the second
for large $\tau$. Since $h_0>0$, initial bodies are ordered by
$\tau$. The comparison principle and
Lemma~\ref{finite time-extinction} show that the extinction set is
downward closed and the expansion set is upward closed. In
particular, a larger solution cannot become extinct while enclosing
a smaller global expanding solution. As the two sets are disjoint,
we obtain
\begin{eqnarray}\label{critical-scales-bounds}
 0<\frac{(\max\limits_{\C_\theta}f)^{-1/n}}
           {\max\limits_{\C_\theta}(h_0/\ell)}
 \leq\tau_{\ast}\leq\tau^{\ast}
 \leq\frac{m^{-1/n}}{\min\limits_{\C_\theta}(h_0/\ell)}<\infty.
\end{eqnarray}

Next, we prove that these sets are open. For an extinguishing
solution, Lemma~\ref{finite time-extinction} supplies a time at which
the solution is enclosed by a strictly subcritical cap. Enlarge this
cap slightly while keeping its radius subcritical. Its support
function then exceeds that of the solution by a positive amount on
the compact set $\C_\theta$. Continuous dependence up to the chosen
time preserves this enclosure for all nearby initial scaling
parameters. The first criterion in \eqref{radius-triggers} proves
openness of the extinction set.

For an expanding solution, choose a time at which it contains a
strictly supercritical cap. Decrease the cap radius slightly while
keeping it supercritical. The resulting positive support function
gap persists for nearby scaling parameters by continuous dependence.
The second criterion in \eqref{radius-triggers} proves openness of
the expansion set. Therefore the two sets in
\eqref{critical-scales} are precisely $(0,\tau_{\ast})$ and
$(\tau^{\ast},\infty)$.

\begin{thm}\label{thm-1}
For every $\tau\in[\tau_{\ast},\tau^{\ast}]$, the solution of \eqref{scalarflow} with
initial support function $\tau h_0$ exists for all time. Moreover, there exists some uniform constant $C>0$ such that 
\begin{eqnarray}\label{non-collapse}
 C^{-1}\leq r(t)\leq R(t)\leq C,
 \qquad t\geq0,
\end{eqnarray}
and
\begin{eqnarray}\label{C2-est}
 \max\limits_{\C_\theta}H(\cdot,t)\leq C,
 \qquad t\geq0.
\end{eqnarray}
\end{thm}

\begin{proof}
For such a parameter $\tau$, neither criterion in
\eqref{radius-triggers} can hold. Thus, throughout the maximal
existence interval,
\begin{eqnarray}\label{uper-low}
 r(t)\leq m^{-1/n},\qquad
 R(t)\geq\left(\max\limits_{\C_\theta}f\right)^{-1/n}.
\end{eqnarray}
By Lemma~\ref{finite time-extinction}, the positive lower bound for
$R(t)$ rules out a finite maximal time. Hence the solution is global.

Fix $T>0$ and choose a time in $[0,T]$ at which $R(t)=R_T$.
Combining \eqref{radius-curvature}, \eqref{C2} and the upper bound
for $r(t)$ at this time, we obtain
\begin{eqnarray*}
 R_T^2\leq C(1+R_T^\gamma).
\end{eqnarray*}
Since $\gamma<2$, this bounds $R_T$ independently of $T$. By
\eqref{C2}, the same is true of $H$. Applying
\eqref{radius-curvature} once more, now with the lower bound for
$R(t)$ in \eqref{uper-low}, yields the assertion \eqref{non-collapse} and \eqref{C2-est}.
 
By \eqref{critical-scales-bounds} and the locally uniform dependence
of the preceding estimates on the initial data, these constants can
be chosen uniformly for $\tau\in[\tau_{\ast},\tau^{\ast}]$.
This completes the proof.
\end{proof}

By Proposition~\ref{Prop-K-upp}, \eqref{non-collapse} also ensures a
time independent upper bound for $K$. Together with \eqref{C2-est},
this yields uniform two-sided bounds for $A$. The regularity argument
in Section~\ref{sec:regularity} then implies uniform higher order
estimates on fixed time intervals after a fixed horizontal
translation. We use these estimates in the convergence argument below.

\section{Convergence after translation}\label{sec:convergence}

Under the uniform estimates established above, the only remaining
loss of compactness is a horizontal translation. We first determine
the translating velocity and then use weighted integrals to control
the position of the translated hypersurfaces.

The following lemma is a weighted spherical cap analogue of \cite[Lemma~4.1]{STW}. We include a proof based on the strict convexity of  \eqref{eta-potential}. 

\begin{lem}\label{lemme-ho1}
Let $g$ be a positive continuous function on $\C_\theta$.
There exists a unique horizontal vector $\eta$ such that
\begin{eqnarray}\label{eta}
 \int_{\C_\theta}\xi_i g(\xi)
 e^{\ell(\xi)^{-1}\<\eta,\xi\>} d\sigma=0,
 \qquad i=1,\cdots,n.
\end{eqnarray}
\end{lem}

\begin{proof}
For a horizontal vector $\eta\in \partial\ov{\RR^{n+1}_{+}}$, we introduce
\begin{eqnarray}\label{eta-potential}
 \Psi(\eta)\coloneqq
 \int_{\C_\theta}\ell(\xi)g(\xi)
 e^{\ell(\xi)^{-1}\<\eta,\xi\>} d\sigma.
\end{eqnarray}
A direct computation yields
\begin{eqnarray}\label{eta-hessian}
 \frac{\partial^2\Psi}{\partial\eta_i\partial\eta_j}
 =\int_{\C_\theta}\frac{\xi_i\xi_j}{\ell(\xi)}g(\xi)
 e^{\ell(\xi)^{-1}\<\eta,\xi\>} d\sigma,
 \qquad 1\leq i,j\leq n.
\end{eqnarray}
%For every nonzero horizontal vector $a$, the function $\<a,\xi\>$ is nonzero on a subset of positive measure in $\C_\theta$. 
Since $\ell$ and $g$ are positive, the Hessian in
\eqref{eta-hessian} is positive definite. Thus $\Psi$ is strictly
convex.

To establish coercivity, write $\eta\coloneqq re_0$, where $r=|\eta|$ and
$e_0\in \SS^{n-1} \cap \p\RR^{n+1}_+$ is a horizontal unit vector. By rotational symmetry, there
exist $\delta,c>0$, independent of $e_0$, such that
$\{\xi\in\C_\theta:\<e_0,\xi\>\geq\delta\}$ has measure at
least $c$. Consequently,
\begin{eqnarray*}
 \Psi(re_0)
 \geq c\min\limits_{\C_\theta}\ell\min\limits_{\C_\theta}g
 \exp\left(\frac{r\delta}{\max\limits_{\C_\theta}\ell}\right)
 \to +\infty\qquad\text{as }r\to +\infty,
\end{eqnarray*}
uniformly in $e_0$. Therefore $\Psi$ has a unique minimizer.
The equation $D\Psi=0$ at this minimizer is precisely Eq. \eqref{eta}.
This completes the proof.
\end{proof}

We are now ready to conclude the proof of Theorem \ref{main-theorem}.
\begin{proof}[\textbf{Proof of Theorem~\ref{main-theorem}}] 
By Theorem~\ref{thm-1}, the solution exists for all time when $\tau\in[\tau_\ast,\tau^\ast]$. We first prove convergence after translation and then show that this interval consists of a single point.

\medskip
\noindent\textit{Step~1. Convergence in the translating frame.} Applying Lemma~\ref{lemme-ho1} to $g=f^{-1}$ and taking the
negative of the resulting vector, we obtain a unique horizontal
vector $\eta$ such that
\begin{eqnarray}\label{ortho}
 \int_{\C_\theta}\frac{\xi_i}{\widetilde f} d\sigma=0,
 \qquad i=1,\cdots,n,\qquad
 \widetilde f\coloneqq f e^{\ell(\xi)^{-1}\<\eta,\xi\>}.
\end{eqnarray}
Fix $\tau\in[\tau_\ast,\tau^\ast]$ and define
$\widetilde X\coloneqq X+t\eta$. Then its support function 
$$\widetilde h\coloneqq h+t\<\eta,\xi\>.$$
We denote the translated hypersurface by $\widetilde\S_t$ and
write $\widetilde A\coloneqq\n^2\widetilde h+
\widetilde h\sigma$. Since horizontal linear functions lie in the
kernel of $\n^2+\sigma$ and satisfy the Robin boundary condition,
we have $\widetilde A=A$ and \eqref{scalarflow} implies 
\begin{eqnarray}\label{translated-flow}
 \begin{cases}
 \partial_t\widetilde h
 =\ell\bigl(\log\det\widetilde A+\log\widetilde f\bigr),
 &\text{in }\C_\theta\times[0, T),\\
 \n_\mu\widetilde h=\cot\theta\widetilde h,
 &\text{on }\partial\C_\theta\times[0,T),\\ 
 \wt h(\cdot,0)=\tau h_0(\cdot),  & \text{on } \C_\theta.
 \end{cases}
\end{eqnarray}
Denote 
$\widetilde F\coloneqq
\log\det\widetilde A+\log\widetilde f
=\ell^{-1}\partial_t\widetilde h$.
By Theorem~\ref{thm-1} and Proposition~\ref{Prop-K-upp},
the eigenvalues of $A$ are uniformly bounded above and below.
Consequently,
\begin{eqnarray}\label{eqn:bound-F-tilde}
|\widetilde F|\leq C. 
\end{eqnarray}

We first establish a uniform bound for the functional associated
with \eqref{translated-flow}.
For each $t$, choose $p_{\text{S}}(t)$ in the relative interior of
$\widehat{\partial\S_t}$ in the supporting hyperplane.
By \eqref{non-collapse}, the diameter of $\widehat\S_t$
is uniformly bounded, and hence
\begin{eqnarray}\label{supp-upbound}
    0<h-\<p_{\text{S}}(t),\xi\>\leq C, \qquad \forall \xi\in \C_\theta. 
\end{eqnarray} 
Define the functional associated with flow \eqref{translated-flow} 
\begin{eqnarray*}
 \widetilde\J(t)
 \coloneqq-\mathrm{Vol}(\widehat{\widetilde\S_t})
 +\int_{\C_\theta}\frac{\widetilde h}{\widetilde f} d\sigma.
\end{eqnarray*}
Since $\eta$ and $p_{\text{S}}(t)$ are horizontal, \eqref{ortho} implies 
\begin{eqnarray*}
 \int_{\C_\theta}
 \frac{\<t\eta+p_{\text{S}}(t),\xi\>}{\widetilde f} d\sigma=0.
\end{eqnarray*}
Using invariance of volume under translations, we therefore obtain
\begin{eqnarray*}
 \widetilde\J(t)
 =-\mathrm{Vol}(\widehat\S_t)
 +\int_{\C_\theta}
 \frac{h-\<p_{\text{S}}(t),\xi\>}{\widetilde f} d\sigma.
\end{eqnarray*}
The diameter bound also controls the enclosed volume; together with \eqref{supp-upbound}, we obtain
\begin{eqnarray}\label{eqn:bound-energy-translated}
     |\widetilde\J(t)|\leq C
\end{eqnarray}where $C$ is independent of $t$.

The flow \eqref{translated-flow} takes the same form
as \eqref{scalarflow}, with a fixed positive prescribed function
$\widetilde f$. Applying \eqref{J} to flow \eqref{translated-flow}, we obtain
\begin{eqnarray}\label{eqn:J-derivative}
 -\widetilde \J'(t)= \int_{\C_\theta}\frac{\ell}{\widetilde f} (e^{\widetilde F}-1)\widetilde F d\sigma 
\end{eqnarray}
By \eqref{eqn:bound-F-tilde}, we have 
\begin{eqnarray*}
 (e^{\widetilde F}-1)\widetilde F
 =\widetilde F^2\int_0^1 e^{s\widetilde F} ds
 \geq C^{-1}\widetilde F^2=C^{-1} \ell^{-2} (\p_t \wt h)^2.
\end{eqnarray*}
Combined with \eqref{eqn:J-derivative} and the positive upper and lower bounds for $\ell$
and $\widetilde f$, this yields
\begin{eqnarray*}
 -\widetilde\J'(t)
 \geq C^{-1}\int_{\C_\theta}
       (\partial_t\widetilde h)^2 d\sigma.
\end{eqnarray*}
In particular, $\widetilde\J$ is nonincreasing. Since it is bounded
below by \eqref{eqn:bound-energy-translated}, integrating in time yields
\begin{eqnarray}\label{C0}
 C  \geq \widetilde\J(0)-\widetilde\J(t)
 =-\int_0^t\widetilde\J'(r)dr \geq C^{-1}\int_0^t\int_{\C_\theta}
       \bigl(\partial_r\widetilde h(\xi,r)\bigr)^2d\sigma dr.
\end{eqnarray}

Next, we prove $\widetilde h(\xi, t)$ is uniformly bounded. We \textbf{claim} that it is equivalent to show
\begin{eqnarray}\label{inte-bound}
    \int_{\C_{\theta}}\frac{\xi_{i}\widetilde h}{\ell \widetilde f}d\s \quad \text{is~bounded~for~}1\leq i\leq n.
\end{eqnarray}
In fact, since
\begin{eqnarray}\label{sup-sum}
\widetilde h  =h-\<p_{\text{S}}(t),\xi\>+\langle p_{\text{S}}(t)+t\eta,\xi\rangle, 
\end{eqnarray}
we have
\begin{eqnarray}\label{eqn:intermide-wt-h}
    \int_{\C_{\theta}}\frac{\xi_{i}\widetilde h}{\ell \widetilde f}d\s=\sum\limits_{j=1}^{n}\<p_{\text{S}}(t)+t\eta, E_{j}\>\int_{\C_{\theta}}\frac{\xi_{i}\xi_{j}}{\ell \widetilde f}d\s+\int_{\C_{\theta}}\frac{\xi_{i}(h-\<p_{\text{S}}(t), \xi\>) }{\ell \widetilde f}d\s,\qquad 1\le i\leq n. 
\end{eqnarray}
The matrix
$\left(\int_{\C_\theta}\frac{\xi_i\xi_j}{\ell\widetilde f}
d\sigma\right)_{1\leq i,j\leq n}$
is positive definite.
Thus, if \eqref{inte-bound} holds, then
\eqref{eqn:intermide-wt-h} and \eqref{supp-upbound} imply
\begin{eqnarray*}
    |p_{\text{S}}(t)+t\eta|\leq C.
\end{eqnarray*}
 In view of \eqref{sup-sum} and using \eqref{supp-upbound} again, we see that \eqref{inte-bound} will imply a uniform $C^0$ bound for $\widetilde h$.  The reverse side is obvious by \eqref{eqn:intermide-wt-h}.

Therefore we proceed to show that \eqref{inte-bound} holds.  From \cite[Proposition 2.6]{MWW-AIM}, we have 
\begin{eqnarray} \label{zero-condi}
    0=\int_{\C_{\theta}}\xi_{i}\det \widetilde A d\s=\int_{\C_{\theta}}\frac{\xi_{i}e^{\widetilde F}}{\widetilde f}d\s, \quad i=1,\cdots, n.
\end{eqnarray}
Combining \eqref{ortho} and \eqref{zero-condi}, we deduce 
\begin{eqnarray}\label{deri}
    \frac{d}{dt}\left(\int_{\C_{\theta}} \frac{\xi_{i}\widetilde h}{\ell \widetilde f}d\s\right)=\int_{\C_{\theta}} \frac{\xi_{i} \widetilde F}{  \widetilde f}d\s=-\int_{\C_\theta}\frac{\xi_i}{\widetilde f}
       (e^{\widetilde F}-1-\widetilde F) d\sigma.
\end{eqnarray}
Together with \eqref{eqn:bound-F-tilde}, we derive 
\begin{eqnarray*}
    0\leq e^{\widetilde F}-1-\wt F =\wt F^2\int_0^1 (1-r)e^{r\wt F} dr \leq  \frac{e^C}2 \wt F^2\leq C_1 (\partial_{t}\widetilde h)^{2}.
\end{eqnarray*}
Integrating  both side of \eqref{deri} over $[0,t]$ and using \eqref{C0},  we obtain
\begin{eqnarray*}
    \left|\int_{\C_{\theta}} \frac{\xi_{i}\widetilde h(\xi,t) }{\ell \widetilde f} d\s\right |\leq \left|\int_{\C_{\theta}}\frac{\xi_{i}\widetilde h(\xi, 0)}{\ell \widetilde f}d\s\right|+    C \int_0^t \int_{\C_{\theta}} \left(\partial_r \widetilde h(\xi,r) \right)^{2}d\s dr\leq C,
\end{eqnarray*}
hence \eqref{inte-bound} holds, and the claim follows, which in turn implies the uniform $C^0$ bound of $\wt h$.

The estimates in Section~\ref{sec:regularity}, applied to
\eqref{translated-flow}, now provide the  uniform higher order bounds for
$\widetilde h$ on $t\geq1$. In particular,
$\{\widetilde h(\cdot,t):t\geq1\}$ is smoothly precompact,
and $\partial_t\widetilde h$ is uniformly continuous in space
and time. It follows from \eqref{C0} that
\begin{eqnarray}\label{speed-decay}
 \|\partial_t\widetilde h(\cdot,t)\|_{C^0(\C_\theta)}
 \to 0
 \qquad\text{as }~~ t\to\infty.
\end{eqnarray}

For any sequence of times tending to infinity, smooth precompactness
provides a subsequence converging to a limit $\widetilde h_\infty$.
Passing to the limit in \eqref{translated-flow} and using
\eqref{speed-decay}, we obtain
\begin{eqnarray}\label{A0}
 \begin{cases}
  \det\widetilde A_\infty=\widetilde f^{-1}, 
 & \qquad \text{in }\C_\theta, \\
  \n_\mu\widetilde h_\infty=\cot\theta\widetilde h_\infty, 
 & \qquad \text{on }\partial\C_\theta,
 \end{cases}
\end{eqnarray}
 where
\begin{eqnarray*}
     \widetilde A_\infty\coloneqq\n^2\widetilde h_\infty+
       \widetilde h_\infty\sigma>0.
\end{eqnarray*}
The strict inequality follows from the uniform lower bound for $A$.
By the uniqueness statement in \cite[Theorem~1.1]{MWW-AIM},
any two such limits differ by a horizontal linear function.

\medskip 

It remains to prove uniqueness of the limit. Since
\begin{eqnarray*}
 0\leq e^x-1-x\leq(e^x-1)x,
 \qquad  \forall ~ x\in\mathbb R,
\end{eqnarray*}
combining \eqref{deri} and \eqref{eqn:J-derivative}, we obtain
\begin{eqnarray*}
 \left|\int_{\C_\theta}
 \frac{\xi_i(\widetilde h(\xi,t)-\widetilde h(\xi,s))}
      {\ell\widetilde f}d\sigma\right| &= &\left|\int_{s}^{t}\frac{d}{dr}\left(\int_{\C_{\theta}} \frac{\xi_{i}\widetilde h(\xi, r)}{\ell \widetilde f}d\s\right) dr\right| \\& \leq &  \int_{s}^{t} \int_{\C_\theta} \frac{|\xi_i|}{\wt f} \left(e^{\wt F}-1-\wt F\right)d\s dr  \\ 
 &\leq& \sup_{\C_\theta} \left(\frac{|\xi_i|}{\ell}\right)  \int_s^t\int_{\C_\theta}
 \frac{\ell}{\widetilde f}
 (e^{\widetilde F}-1)\widetilde Fd\sigma dr\\
 &=&C \bigl(\widetilde\J(s)-\widetilde\J(t)\bigr),
 \qquad t\geq s.
\end{eqnarray*}
By \eqref{eqn:bound-energy-translated} and
\eqref{eqn:J-derivative}, $\widetilde\J$ has a finite limit.
The right-hand side therefore tends to zero as $s,t\to\infty$.
Consequently, any two subsequential limits
$\widetilde h_\infty$ and $\widetilde h^\infty$ satisfy
\begin{eqnarray*}
 \int_{\C_\theta}
 \frac{\xi_i(\widetilde h_\infty-\widetilde h^\infty)}
      {\ell\widetilde f}d\sigma=0,
 \qquad i=1,\ldots,n.
\end{eqnarray*}
 with
\begin{eqnarray*}
     \widetilde h_\infty-\widetilde h^\infty
=\sum_{j=1}^n a_j\xi_j.
\end{eqnarray*}
Consequently,
\begin{eqnarray*}
 0=\sum_{i=1}^n a_i\int_{\C_\theta}
 \frac{\xi_i(\widetilde h_\infty-\widetilde h^\infty)}{\ell\widetilde f} d\sigma=\int_{\C_\theta}
 \frac{\bigl(\sum_{i=1}^n a_i\xi_i\bigr)^2}
      {\ell\widetilde f} d\sigma.
\end{eqnarray*}
Since $\left(\int_{\C_{\theta}}\frac{\xi_{i}\xi_{j}}{\ell \widetilde f}d\s\right)_{1\leq i,j\leq n} $ is positive definite, all $a_i$ must vanish.
Together with the $C^\infty$ precompactness, this yields
$\widetilde h(\cdot,t)\to\widetilde h_\infty$ smoothly. Hence
\begin{eqnarray*}
 \widetilde X(\xi,t)
=\n\widetilde h(\xi,t)
    +\widetilde h(\xi,t)(\xi-\cos\theta e)\to 
 \n\widetilde h_\infty(\xi)
 +\widetilde h_\infty(\xi)(\xi-\cos\theta e)
 \eqqcolon X^* \quad\hbox{in }C^\infty(\C_\theta),
\end{eqnarray*}
which proves \eqref{asy}. By \eqref{A0}, the limiting hypersurface $X^*$ has Gauss curvature
\begin{eqnarray*}
 K_\infty=\bigl[\det(\n^2\widetilde h_\infty
              +\widetilde h_\infty\sigma)\bigr]^{-1}=\widetilde f
 =f e^{\ell(\xi) ^{-1}\<\eta,\xi\>}.
\end{eqnarray*}

\medskip

\noindent\textit{Step~2. Uniqueness of the critical scaling.}
Let $h^\ast(\xi,t)$ and $h_\ast(\xi,t)$ be the solutions starting
from $\tau^\ast h_0$ and $\tau_\ast h_0$, respectively.
By \textit{Step~1}, their translated limits are
\begin{eqnarray*}%\label{endpoint-limits}
  h^\ast(\xi)
 \coloneqq \lim_{t\to\infty}
 \left(h^\ast(\xi,t)+t\<\eta,\xi\>\right),\qquad  
 h_\ast(\xi)
 \coloneqq \lim_{t\to\infty}
 \left(h_\ast(\xi,t)+t\<\eta,\xi\>\right).
\end{eqnarray*}
Both solve \eqref{A0} with the same $\widetilde f$, and hence
$h^\ast-h_\ast=\sum_{i=1}^n a_i\xi_i$.
By the normalization in Section~\ref{sec:selection}, the origin
lies in the relative interior of the initial planar wetting part. Thus
$\varepsilon_0E_{n+1}\in\widehat \S_0$ for some $\varepsilon_0>0$, and 
\begin{eqnarray}\label{eqn:normalized-positive-h-0}
 h_0\geq\varepsilon_0\cos\theta>0,
 \qquad
\text{ and } ~~ \min_{\C_\theta}\frac{h_0}{\ell}>0.
\end{eqnarray}
Since $\tau_\ast\leq\tau^\ast$, Lemma~\ref{lem-maxi} yields 
$h^\ast(\xi,t)\geq h_\ast(\xi,t)$ and thus
$h^\ast\geq h_\ast$. A horizontal linear function cannot be
nonnegative on the whole rotationally symmetric cap unless it
vanishes. Therefore $h^\ast=h_\ast$.

Denote 
$$\omega(\cdot, t) \coloneqq(h^\ast(\cdot,t)-h_\ast(\cdot,t))/\ell,$$ then the proof argument in Lemma~\ref{lem-maxi} yields
\begin{eqnarray}\label{endpoint-difference}
 \begin{cases}
\partial_t\omega=a^{ij}\omega_{ij}+b^i\omega_i+c\omega,
 & \qquad \text{in }\C_\theta\times(0,\infty),\\
\n_\mu\omega=0,
 & \qquad \text{on }\partial\C_\theta\times(0,\infty),\\ 
 \omega(\cdot,0)=(\tau^*-\tau_*)h_0/\ell,  & \qquad  \text{on } \C_\theta.
 \end{cases}
\end{eqnarray}
On every finite time interval, the constant function 
$C\coloneqq(\tau^\ast-\tau_\ast)
\min\limits_{\C_\theta}(h_0/\ell)\geq0$ is a subsolution of
\eqref{endpoint-difference}. Applying the comparison principle implies
\begin{eqnarray*}%\label{endpoint-separation}
 h^\ast(\xi,t)-h_\ast(\xi,t)
 \geq\ell(\xi)(\tau^\ast-\tau_\ast)
 \min\limits_{\C_\theta}\frac{h_0}{\ell},
 \qquad t\geq0.
\end{eqnarray*}
Passing to the limit and using $h^\ast=h_\ast$, we find
\begin{eqnarray*}%\label{threshold-unique}
 0\geq\ell(\xi)(\tau^\ast-\tau_\ast)
 \min\limits_{\C_\theta}\frac{h_0}{\ell}.
\end{eqnarray*}
Together with \eqref{eqn:normalized-positive-h-0}, this implies
$\tau_\ast=\tau^\ast$. Hence the proof is complete.
\end{proof} 
 
\printbibliography

\end{document}